\documentclass[reqno,12pt,a4paper]{amsart}
\usepackage{amsmath,amsthm,amssymb,tabularx}
\usepackage[cmtip,arrow,all,2cell]{xy}

\theoremstyle{plain}
\newtheorem{tm}{Theorem}[section]
\newtheorem{lm}[tm]{Lemma}
\newtheorem{cor}[tm]{Corollary}
\newtheorem{prop}[tm]{Proposition}

\theoremstyle{definition}
\newtheorem{definition}[tm]{Definition}
\newtheorem{ex}[tm]{Example}
\newtheorem{rem}[tm]{Remark}
\newtheorem{nots}[tm]{Notations}

\newcommand{\beq}{\begin{equation}}
\newcommand{\eeq}{\end{equation}}
\newcommand{\bga}{\begin{gather*}}
\newcommand{\ega}{\end{gather*}}
\newcommand{\bal}{\begin{align*}}
\newcommand{\eal}{\end{align*}}
\newcommand{\bit}{\begin{itemize}}
\newcommand{\eit}{\end{itemize}}
\newcommand{\btm}{\begin{tm}}
\newcommand{\etm}{\end{tm}}
\newcommand{\blm}{\begin{lm}}
\newcommand{\elm}{\end{lm}}
\newcommand{\bcor}{\begin{cor}}
\newcommand{\ecor}{\end{cor}}
\newcommand{\bprop}{\begin{prop}}
\newcommand{\eprop}{\end{prop}}
\newcommand{\bex}{\begin{ex}}
\newcommand{\eex}{\end{ex}}
\newcommand{\bpr}{\begin{proof}}
\newcommand{\epr}{\end{proof}}
\newcommand{\brem}{\begin{rem}}
\newcommand{\erem}{\end{rem}}

\def\C{\mathbb{C}}
\def\N{\mathbb{N}}
\def\R{\mathbb{R}}

\def\Id{\mathbb I}

\def\K{{\mathcal K}}

\def\e{\varepsilon}
\let\a\alpha
\def\trans#1#2{{}^{#1}\mkern-3mu #2}
\let\cal\mathcal
\def\m{\mathfrak m}
\def\hatC{\widehat C}
\def\id{{\rm id}}

\def \le {\leqslant}
\def \ge {\geqslant}
\let \<\langle
\let \>\rangle
\let\phi\varphi
\let\kappa\varkappa
\def\hathatA{\skew4\widehat{\widehat A}}

\def\cst{$C^*$}
\def\omax{\underset\max\otimes}
\def\omin{\underset\min\otimes}
\def\liminv{\mathop{\underset\longleftarrow\lim}}

\def\supp{\operatorname{supp}}
\begin{document}
\author{Yulia Kuznetsova}
\thanks{Supported by the research grant F1R-MTH-PUL-10NCHA of the University of Luxembourg}
\title{A duality for Moore groups}
\address{University of Luxembourg
6, rue Richard Coudenhove-Kalergi
L-1359 Luxembourg
Tel: (+352) 466644 6802}
\email{julia.kuznetsova@uni.lu}
\maketitle

\begin{abstract}
We suggest a new generalization of Pontryagin duality from the category of Abelian locally compact groups to a category which includes all Moore groups, i.e. groups whose irreducible representations are finite-dimensional. Objects in this category are pro-C*-algebras with a structure of Hopf algebras (in the strict algebraic sense) with respect to a certain topological tensor product.
\end{abstract}

\section{Introduction}

This paper contributes to the problem of duality of locally compact groups. For history and background of this question, we refer to monographs \cite{enock} and \cite{hajac}.

It is known that in the theories of Kac algebras and of locally compact quantum groups, one has to require axiomatically the existence of invariant weights, which are analogues of the Haar measure of a group. This need is due to the very definition of the duality functor, which utilizes a distinguished {\it regular} representation of a Kac algebra or of a locally compact quantum group.

It would be desirable to have a category with a duality functor that does not involve the Haar measure/weight. While this task is quite far from the solution in the general case, in particular cases one can construct dualities which behave very well.

The basic example is of course the category of Abelian locally compact groups and the Pontryagin duality for them. Another simple case is that of finite groups and their function algebras with pointwise multiplication and with convolution, regarded as finite-dimensional Hopf algebras.

In this paper we suggest a ``continuous non-Abelian'' example of a duality functor, independent of the Haar measure. It acts on a category which includes all Abelian and compact locally compact groups, and preserves the Pontryagin duality for the Abelian case. This approach is inspired by a recent work by Akbarov \cite{akb-hopf} who constructs a duality for complex Lie groups with algebraic component of identity.

The objects of our category are pro-\cst-algebras (called also locally \cst-algebras), that is, inverse, or projective, limits of \cst-algebras. For any topological algebra $A$ with involution, there is a universal pro-\cst-algebra containing $A$; it is called the \cst-envelope $A^\Diamond$ of $A$ (see section \ref{section-envelopes}).

We consider a subclass of pro-\cst-algebras which have a structure of a Hopf algebra with involution with respect to a topological tensor product $\odot$. For rather wide class of topological vector spaces, this is the same as injective tensor product, see further details in section \ref{section-Moore}. These Hopf algebras are called $\odot$-Hopf *-algebras.

For a $\odot$-Hopf *-algebra $A$, the dual algebra is defined as follows. The space $A^\star$ of linear continuous functionals on $A$, taken with the topology of uniform convergence on totally bounded subsets of $A$, has a natural algebra structure, so we can define its \cst-envelope $(A^\star)^\Diamond$. This is by definition the dual algebra $\widehat A$ of $A$. The second dual algebra is obtained by repeated application of this procedure: $\hathatA=\big(({\widehat A})^\star\big)^\Diamond$, and $A$ is said to be reflexive if $\hathatA=A$.

The main result of this paper is proved in section~\ref{section_SIN}: there exists a category of reflexive $\odot$-Hopf *-algebras, which contains the algebras $C(G)$ for all Moore groups $G$. On this category, $\widehat{}\;$ is a duality functor.

We would like to discuss the group case in more detail. Recall that $G$ is called a Moore group if all its irreducible representations are finite-dimensional. This class includes all Abelian and compact groups but not all discrete groups. The diagram \eqref{quadrat-intro} below illustrates the duality construction in the case of a Moore group.

We start with the algebra $A=C(G)$ of all continuous functions on $G$, in compact-open topology (note that this is not a Banach algebra unless $G$ is compact). Then we pass to its conjugate space $C(G)^\star$, which is the space $M_c(G)$ of all compactly supported measures. The \cst-envelope $\hatC(G)$ of $M_c(G)$ is by definition the dual algebra of $C(G)$. In its turn, the conjugate space $\hatC(G)^\star=\mathcal K(G)$ may be identified with a subalgebra in $C(G)$. In sections \ref{section-K(G)} and \ref{section-stereo} we study properties of $\mathcal K(G)$ and $\hatC(G)$, respectively. In theorem \ref{SIN-reflexive} we prove that the \cst-envelope of $\K(G)$ is again $C(G)$, what is illustrated by the following diagram:
\beq\label{quadrat-intro}
 \xymatrix @R=1.pc @C=2.pc
 {
 C(G)\ar@{|->}[rrr]^{\star}
 &
 & &
 M_c(G) \ar@{->}[dd]^{C^*\rm -env}
 \\
 & & &
 \\
 \K(G) \ar@{->}[uu]^{C^*\rm -env}
 & & &
 \hatC(G)\ar@{|->}[lll]_{\star}
 }
\eeq

\section{\cst-envelopes}\label{section-envelopes}

Our main objects are pro-\cst-algebras, or locally \cst-algebras, as they are also called. The theory of these algebras was developed by J.~Inoue \cite{inoue}, M.~Fragoulopoulou \cite{fragoul-envel}, C.~Phillips \cite{phil} and other authors.

\begin{definition}
A {\it \cst-seminorm } on an algebra $A$ with involution is such a seminorm $p$ that $p(x^*x)=p(x)^2$ for all $x\in A$. We say also that $p\,$ possesses the {\it \cst-property}.
\end{definition}

This definition says nothing about continuity of multiplication and involution, but due to a theorem of Sebasty\'en \cite{seb} every \cst-seminorm $p$ automatically satisfies $p(xy)\le p(x)p(y)$ and $p(x^*)=p(x)$ for any $x,y\in A$. This means that multiplication and involution are automatically continuous (with respect to $p$).

\begin{definition}
A {\it pro-\cst-algebra\/} is a complete Hausdorff topological algebra $A$ with jointly continuous multiplication with involution such that its topology is generated by a family of \cst-seminorms.
\end{definition}

Every pro-\cst-algebra $A$ is equal to the inverse, or projective, limit of (Banach) \cst-algebras $A_p$ over all $p$, where $A_p$ is obtained as the $p$-completion of the quotient algebra $A/\ker p$ \cite[Proposition 1.2]{phil}.

Next we will give a definition of the \cst-envelope of a topological algebra, which is the biggest pro-\cst-algebra containing a given one as a dense subalgebra.
Other authors \cite{inoue,fragoul-envel} define \cst-envelopes only for topological *-algebras with sub\-mul\-ti\-pli\-ca\-tive seminorms (Arens-Michael algebras). This is done to have a correspondence between representations of such an algebra $A$ and its \cst-envelope. But we do not indent to study the representations of $A$, and therefore it is not necessary to impose any additional conditions on it. We require that $A$ is unital, but the multiplication and involution, in principle, can be even discontinuous in the initial topology.
The following notations will be kept throughout the paper.

\begin{nots}\label{nots-envelope}
Let $A$ be a unital algebra with involution, endowed with some locally convex topology.
Denote by $\cal P(A)$ the set of all continuous $C^*$-seminorms on $A$. For every $p\in \cal P(A)$ the kernel $\ker p$ is a *-ideal in $A$, so $A/\ker p$ is a normed algebra with the quotient norm $\bar p$ of $p$. Its completion with respect to $\bar p$ will be denoted by $C_p^*(A)$, or simply $A_p$ when there can be no confusion. This is a (Banach) $C^*$-algebra. The canonical mapping $A\to C_p^*(A)$ will be denoted by $i_p$.
\end{nots}

\blm
The algebras $A_p$ with pointwise ordering on $\mathcal P(A)$ form an inverse spectrum of\/ \cst-algebras.
\elm

\bpr
First of all, $\cal P(A)$ is directed by the pointwise ordering $\le$: for $p,q\in\cal P(A)$ one may put $r=\max(p,q)$, then $p\le r$, $q\le r$ and \cst-property holds for $r$.

If $p,q\in\cal P(A)$ and $p\le q$, then the map $\pi_p^q:i_q(A)\to i_p(A)$, $i_q(x)\mapsto i_p(x)$, is continuous on $i_q(A)$ and can be thus extended to $A_q$. Denote this extension still $\pi_p^q$.

Let now $p\le q\le r$. Obviously on $i_r(A)$ the maps $\pi_p^r$ and $\pi_p^q\circ\pi_q^r$ coincide and map $i_r(x)$ to $i_p(x)$. As $i_r(A)$ is dense in $A_r$, these maps must also coincide on the entire algebra $A_r$. Thus, $\{A_p:p\in\mathcal P(A)\}$ with the maps $\pi^q_p$ form an inverse system.
\epr

\begin{definition}
The {\it $C^*$-envelope} $A^\Diamond$ of an algebra $A$ is the inverse limit of the algebras $A_p$
 over $p\in \mathcal P(A)$.
\end{definition}

The \cst-envelope can be alternatively defined as the completion of $A/E$ with respect to all continuous \cst-seminorms, where $E$ is the common kernel of $p\in \mathcal P(A)$.
Any algebra $A$ is continuously, but not always injectively mapped into its envelope $A^\Diamond$. It may happen that there are no other \cst-seminorms except zero; then $A^\Diamond=\{0\}$.
Note that for the unit 1 of $A$, $p(1)^2=p(1^*1)=p(1)$, so that $p(1)$ must be either 0 or 1. But if it is 0, then $p$ is identically zero, so $p(1)=1$ for every nontrivial seminorm. Note also that if $B\subset A$ is a dense *-subalgebra (with the induced topology), then $B^\Diamond = A^\Diamond$.

\brem\label{universality}
The \cst-envelope has the following universal property. If $\phi: A\to B$ is a continuous *-homomorphism to a \cst-algebra $B$, then it may be uniquely extended to a continuous *-homomorphism $\bar \phi: A^\Diamond\to B$. From this it follows that if $\phi: A\to B$ is a con\-ti\-nu\-ous *-homomorphism of topological algebras, it may be extended to a continuous *-homomorphism $\bar \phi: A^\Diamond\to B^\Diamond$.
\erem

When working with groups, we will use two main pro-\cst-algebras. First of them is $C(G)$, the algebra of all continuous functions on $G$ with the topology of uniform convergence on compact sets. The second one is its dual algebra $\hatC(G)$, defined below as the \cst-envelope of another group algebra. Before giving this definition, we need to discuss the theory of so called stereotype spaces.

\begin{nots}
If $X$ is a locally convex space, $X^\star$ denotes the space of all linear continuous functionals on $X$ in the topology of uniform convergence on totally bounded subsets of $X$. If $X$ is isomorphic to $(X^\star)^\star$, $X$ is called a {\it stereotype space}.
\end{nots}

We use the term ``conjugate space'' for the space of linear functionals to distinguish this notion from the dual algebra introduced in this paper. So, we always consider conjugate spaces in special topology---of uniform convergence on totally bounded sets. It turns out that then almost all classical spaces become reflexive: $X\simeq (X^\star)^\star$; in particular, every Banach or Fr\'echet space, or the space $C(G)$ for a locally compact group $G$. Important results in this direction were obtained already in 1950s by M.~Smith \cite{smith} and 1970s by K.~Brauner \cite{brauner}, but an extensive theory including applications to topological algebra is quite recent and due to S.~Akbarov \cite{akb, akb-hopf}.

Thus, we consider $M_c(G)=C(G)^\star$, the space of all compactly supported measures, with the topology of uniform convergence on totally bounded subsets of $C(G)$. It may be noted that if $X$ is known to be stereotype, then this topology on $X^\star$ coincides with that of uniform convergence on compact sets. This is true, in particular, for $C(G)$ and $M_c(G)$. We will use the following properties of $M_c(G)$ (\cite{akb}, Theorem 10.11 and Example 10.7):

\bprop\label{G-imbedded}
Let $G$ be a locally compact group. Then:\\
(i) $G$ is homeomorphically imbedded into $M_c(G)$ via delta-functions, $t\mapsto \delta_t$;\\
(ii) the linear span of delta-functions is dense in $M_c(G)$.
\eprop

\begin{nots}
The \cst-envelope $M_c(G)^\Diamond$ of $M_c(G)$ will be denoted by $\hatC(G)$. The set $\mathcal P(M_c(G))$ will be also denoted by $\mathcal P(G)$. For $p\in\mathcal P(G)$ we denote $C^*_p(G)=C^*_p(M_c(G))$, and also $C^*_\pi(G)$ is $p$ is the norm of a representation $\pi$. Thus, $\hatC(G)=\liminv C^*_p(G)$. Note that we do not deal with the reduced \cst-algebra of $G$, and the notation $C^*_r(G)$, if used, means just the \cst-algebra associated to a seminorm $r\in\mathcal P(G)$. If $G$ is an Abelian or compact locally compact group, $\widehat G$ denotes its dual group or dual space respectively.
\end{nots}

Every \cst-seminorm $p\in \cal P(G)$ defines a *-homomorphism of $M_c(G)$ into the \cst-algebra $C^*_p(G)$, which may be further mapped into the algebra $\mathcal B(H)$ of bounded operators on a Hilbert space $H$. Let us denote this homomorphism $\psi: M_c(G)\to \mathcal B(H)$. The unit $\delta_e$ of $M_c(G)$ is mapped into a projection, but one can assume that it is mapped into $1_H$, reducing $H$ if necessary. This reduction will not change $p$. Now, as $\delta_t*\delta_t^*=\delta_t^* *\delta_t=\delta_e$, we have $\psi(\delta_t)\psi(\delta_t)^*=\psi(\delta_t)^*\psi(\delta_t)=1_H$, i.e. all operators $\Psi_t=\psi(\delta_t)$ are unitary. Thus, $p$ generates a unitary representation of $G$, and it is important to note that $\Psi$ is norm continuous.

Conversely, every norm continuous representation of $G$ generates (via delta-functions) a non-degenerate norm continuous representation of $M_c(G)$, and consequently, a \cst-seminorm on $M_c(G)$ (see \cite[10.12]{akb}).

Two representations of $G$ may not be equivalent even if they generate the same seminorm, for example, $\pi$ and $\pi\oplus\pi$.

\blm\label{hom-of-hatC(G)}
Let $\phi: G\to H$ be a continuous homomorphism of locally compact groups. Then there is a continuous *-homomorphism $\widehat\phi: \hatC(G)\to\hatC(H)$ such that $\widehat\phi(\delta_t)=\delta_{\phi(t)}$.
\elm
\bpr
Consider the canonical embedding $i:H\to M_c(H)$, $i(t)=\delta_t$. Then $i\circ\phi$ is a continuous homomorphism of $G$ into a stereotype algebra $M_c(H)$, so it is extended to $M_c(G)$ by linearity and continuity \cite[Theorem~10.12]{akb}. Clearly this morphism is involutive. By universality property \ref{universality}, $i\circ\phi$ is extended to a continuous *-morphism of the envelopes $\widehat\phi:\hatC(G)\to\hatC(H)$.
\epr

Now we will describe $\hatC(G)$ in the particular cases of Abelian and compact groups.

\btm\label{C(G)-abel}
Let $G$ be an Abelian locally compact group, and let $\widehat G$ be its dual group.
Then $\hatC(G)$ coincides with the algebra $C(\widehat G)$, taken with the topology of uniform convergence on compact sets.
\etm

\bpr
Denote for brevity $A=M_c(G)$. It is evident that for every compact set $K\subset \widehat G$ the seminorm
$$
p_K(\mu)=\max_{t\in K}|\hat\mu(t)|,
$$
where $\hat\mu$ is the Fourier transform of $\mu$, has the  $C^*$-property. We will show that every other $C^*$-seminorm $p$ equals $p_K$ for some compact $K\subset \widehat G$.

For every $p$, $A_p$ is a unital $C^*$-algebra. As it is commutative, it is isomorphic to the function algebra $C(\Omega)$ on a compact space $\Omega$.  Every point $\omega\in\Omega$ defines a continuous character of the algebra $A_p$; and since the map $A\to A_p$ is continuous, the restriction of $\omega$ onto $A$ is also a continuous character of $A$.

All such characters are just evaluations at the points of the dual group $\widehat G$, and this correspondence with $\widehat G$ is a homeomorphism \cite[10.12]{akb}. Thus, we get an imbedding $\psi:\Omega\to\widehat G$, with
$$p(\mu)=\max_{\omega\in\Omega}|\mu(\omega)|=\max_{t\in\psi(\Omega)} |\hat\mu(t)|,$$
i.e. $p=p_{\psi(\Omega)}$. This proves the theorem.
\epr

To describe $\hatC(G)$ in the compact case, we need the following \cite[Corollary 2]{shtern}:

\blm\label{compact-finite}
Let $T$ be a norm continuous unitary representation of a compact group $G$. It may be decomposed into the direct sum of irreducible unitary representations of \,$G$, where only a finite number of representations are non-equivalent.
\elm

\btm\label{direct_product_irr}
For a compact group $G$, the algebra $\hatC(G)$ is equal to the direct product of the matrix algebras $C^*_\pi(G)$ over all irreducible representations $\pi$.
\etm
\bpr
By lemma \ref{compact-finite}, for every $\pi$ there is a finite set $s(\pi)\subset\widehat G$ such that $\pi=\oplus_{j\in s(\pi)} \pi_j$, and hence $C^*_\pi(G)=\prod_{j\in s(\pi)} C_{\pi_j}^*(G)$. Clearly, $\|\pi\|=\max_{j\in s(\pi)}\|\pi_j\|$. Let us show that if $p=\|\pi\|$, $q=\|\tau\|$ and $p\le q$, then $s(\pi)\subset s(\tau)$. This is sufficient to prove in the case when $\pi$ is irreducible, i.e. $s(\pi)=\{\pi\}$. Since $p\le q$, we have an epimorphism $R: C_\tau^*(G)\to C^*_\pi(G)$. Here $C^*_\pi(G)$ and all $C_{\tau_j}^*(G)$, $j\in s(\tau)$, are just full matrix algebras. The image of every ideal $I_j=C_{\tau_j}^*(G)\times\prod_{k\ne j} \{0\}$ is either $\{0\}$ or $C_\pi^*(G)$, since $C^*_\pi(G)$ is a simple algebra. If there are distinct $j\ne k$ such that $R(I_j)=R(I_k)=C_\pi^*(G)$, then we can take $a\in I_j$, $b\in I_k$ such that $R(a)=R(b)=\Id$; then $R(a+b)=2\Id$, but $R(a+b)=R(a+b)R(a)=R((a+b)a)=R(a)=\Id$. This contradiction shows that such $j$ that $R(I_j)=C^*_\pi(G)$ is unique. The restriction of $R$ to $I_j$ is an isomorphism, since $I_j$ is also a simple algebra. Thus, $\pi$ is equivalent to $\tau_j$, and $s(\pi)=\{\pi\}\subset s(\tau)$, as claimed.

From the other side, for every finite set $s\subset\widehat G$ there is a norm continuous representation $\pi=\oplus_{j\in s}\pi_j$ such that $s(\pi)=s$.

The preceding reasoning shows that $p\mapsto s(\pi)$, where $\pi$ is any representation such that $p=\|\pi\|$, is an order-preserving bijection of $\mathcal P(G)$ onto the set of finite subsets of $\widehat G$, ordered by inclusion. Thus, $\hatC(G)=\liminv_{p\in\mathcal P(G)}C^*_p(G)$ is equal to the inverse limit $\liminv_{s\subset\widehat G} \prod_{\pi\in s} C^*_\pi(G)$. Clearly the latter limit is just the direct product of all $C^*_\pi(G)$, and we have finally $\hatC(G)=\prod_{\pi\in \widehat G}C^*_\pi(G)$.
\epr

\section{The space of coefficients and SIN-groups}\label{section-K(G)}

Now we will discuss the conjugate space $\K(G)=\big( \hatC(G) \big)^\star$ of $\hatC(G)$. Since $G$ is continuously imbedded into $M_c(G)$, $t\mapsto\delta_t$, every functional $f\in \K(G)$ restricts to a continuous function on $G$. Moreover, if two functionals are equal on $\delta$-functions, then by linearity and density they are equal also on $M_c(G)$ and further on the entire algebra $\hatC(G)$. We identify therefore $f$ with the function defined in this way, and represent $\K(G)$ as a subalgebra in $C(G)$. As usual, we endow $\K(G)$ with the topology of uniform convergence on totally bounded sets in $\hatC(G)$.

Recall that an inverse limit $X=\liminv X_\a$ is called {\it reduced\/} if the canonical projection of $X$ is dense in every $X_\a$. By definition the inverse limit $A=\liminv A_p$, $p\in\mathcal P(A)$, is reduced. We will need the following simple result, close to a classical one \cite[Theorem~4.4]{schaefer}:

\bprop\label{conjugate_ind_limit}
Let $X=\liminv X_\a$ be a reduced inverse limit of locally convex spaces. Then, as a set, $X^\star=\cup X_\a^\star$, and all inclusions $X^\star_\a\hookrightarrow X^\star$ are continuous. If a set $B\subset X^\star$ is dense in every $X_\a^\star$ then it is dense in $X^\star$.
\eprop
\bpr
The equality $X^\star=\cup X_\a^\star$ is classical \cite[Theorem~4.4]{schaefer}. The inclusion $X^\star_\a\hookrightarrow X^\star$ is dual to the canonical projection $\pi_\a:X\to X_\a$, and so it is continuous.

It follows that the topology of $X^\star$ is not stronger than the inductive limit topology. It is easy to see that the set $B$ as in the statement is dense in $X$ in the inductive limit topology; as a consequence, it is dense in the topology of $X^\star$ as well.
\epr

\bprop\label{Astar-properties}
Let $A=\liminv\limits_{p\in\mathcal P(A)} A_p$ be a pro-\cst-algebra. Then:

\vskip-14pt
\begin{enumerate}
\item $A^\star=\cup A_p^\star$;
\item if a set $B\subset A^\star$ is dense in every $A_p^\star$ then it is dense in $A^\star$;
\item $A^\star$ is the linear span of the states of all $A_p$;
\item the linear span of pure states of $A_p$ is dense in $A^\star$.
\end{enumerate}
\eprop
\bpr
(1) and (2) are just repetitions of proposition \ref{conjugate_ind_limit}. (3) follows from the classical fact that every continuous linear functional on a \cst-algebra is a finite linear combination of states \cite[4.3.7]{kad}.

(4): It is known that the linear span of pure states on $A_p$ is $A_p$-weakly dense in $A_p^\star$. The unit ball $B_p$ of $A_p^\star$ is compact in the topology of uniform convergence on totally bounded sets in $A_p$ \cite[Example 4.6]{akb}, and this topology is stronger than the weak topology in which $B_p$ is compact too. It follows that these two topologies on $B_p$ coincide. Thus, the linear span of pure states is dense in every $B_p$, and therefore in every $A_p^\star$. By lemma \ref{conjugate_ind_limit} it is also dense in $A^\star$.
\epr

In application to the group case, we get the following:

\bprop\label{K(G)-properties}
$\K(G)$ has the following properties:
\begin{enumerate}
\item $\K(G)$ is contained in the linear span of all positive-definite functions on $G$;
\item $\K(G)$ is equal to the space of the coefficients of all norm-continuous unitary representations of $G$;
\item $\K(G)$ contains the constant 1;
\item all functions $f\in\K(G)$ are bounded;
\item the linear span of the coefficients of irreducible norm continuous representations is dense in $\K(G)$.
\end{enumerate}
\eprop
\bpr
By proposition \ref{Astar-properties}(3), $\K(G)$ is the linear span of states of all $C_p^*(G)$. For every state $f\in \big(C_p^*(G)\big)^\star$, the corresponding function on $G$ is positive definite. This implies (1) and (this is a known fact) (4). From the other side, $f$ has form $f(\mu)=\langle \pi(\mu)x,x\rangle$, $\mu\in C_p^*(G)$, for some representation $\pi$ of $C_p^*(G)$, and for some vector $x$. Then $f$ is a coefficient of the representation $\pi\circ i_p$ of $G$ (see notations \ref{nots-envelope}), if we identify $G$ with its image in $M_c(G)$. Since $\|\pi\circ i_p\|\le p$, this representation is norm continuous and this proves (2). For (3), note that constant 1 function is the coefficient of the trivial representation, which is norm continuous. And finally, (5) follows from proposition \ref{Astar-properties}(4), since coefficients of irreducible representations are pure states.
\epr

\bex
{\bf Discrete groups.} If $G$ is discrete, algebras $C_p^*(G)$ correspond to {\it all\/} representations of $G$, as they are all continuous. This means that $\K(G)$ is the space of the coefficients of all unitary representations of $G$, i.e., the Fourier-Stieltjes algebra $B(G)$. In particular, consider the regular representation on the space $L_2(G)$: $U_s\phi(t) = \phi(s^{-1}t)$. Take $\phi=I_{\{e\}}$, $\psi=I_{\{s\}}$ (indicator functions of the points $e$ and $s$), then the function
$$f(t)=\langle U_t\phi,\psi\rangle = \langle I_{\{t\}},I_{\{s\}}\rangle = I_{\{s\}}(t)$$
belongs to $\K(G)$. Thus, $\K(G)$ contains all finitely supported functions on $G$. On a non-discrete group, the regular representation is not norm continuous and $\K(G)$ does not contain all its coefficients.
\eex

\blm\label{K(G)-algebra}
$\K(G)$ is closed under pointwise multiplication and complex conjugation.
\elm
\bpr
Let us denote $\K_p(G)=\big(C_p^*(G)\big)^\star$, then it follows that $\K(G)=\cup_p\, \K_p(G)$.
Take $f,g\in \K(G)$ and prove that $f\cdot g\in \K(G)$. There exist \cst-seminorms $p,q$ such that $f\in \mathcal K_p(G)$, $g\in \mathcal K_q(G)$. Consider any \cst-tensor product $\cal C=C_p^*(G)\otimes C_q^*(G)$ of the algebras $C_p^*(G)$ and $C_q^*(G)$. The mapping $U:G\to \cal C$, $t\mapsto i_p(\delta_t)\otimes i_q(\delta_t)$ (see notations \ref{nots-envelope}), is a unitary representation of $G$ in the Hilbert space where $\cal C$ acts non-degenerately. It generates a representation of $M_c(G)$ and thus a continuous \cst-seminorm $r$ on $M_c(G)$, so that $\mathcal C=C_r^*(G)$.

By proposition \ref{Astar-properties}(3), we can assume that $f$ and $g$ are states. Then there exists \cite[11.3]{kad} a product state  $\tau$ on $\mathcal C$ such that $\tau(x\otimes y)=f(x)\cdot g(y)$ for all $x\in C_p^*(G)$, $y\in C_q^*(G)$. For $t\in G$ we have then $\tau(t)=\tau(i_r(\delta_t))=f(t)\cdot g(t)$, so $fg\in\K(G)$ and we have proved the first statement.

The second statement is proved exactly as in general theory (e.g. \cite[27.26]{HR}).
\epr

Our aim is the commutative duality diagram \eqref{quadrat-intro}, and for this we need that $\K(G)^\Diamond=C(G)$. To have this equality, it is necessary (but not sufficient) that $\K(G)$ separates points of $G$. A convenient class of groups where this is true is the class of SIN-groups.

These are groups with a basis of invariant neighborhoods of identity, i.e. such neighbourhoods $U$ that $g\,Ug^{-1}=U$ for all $g\in G$. These groups were introduced and studied in detail by Grosser and Moskowitz \cite{gros-mosk}.

The class of SIN-groups includes all Abelian, compact and discrete groups. There is an equivalent definition: a locally compact group is a SIN-group if and only if its left and right uniform structures are equivalent. In particular, such groups are always unimodular. For connected groups, this class coincides with that of MAP-groups (maximally almost periodic, or groups on which finite-dimensional unitary representations separate points) and Z-groups (such that the quotient group over the center is compact). Moreover, a connected SIN-group is a direct product of $\R^n$ by a compact group. In general, there is a complete structure theorem for SIN-groups \cite{gros-mosk}: a group $G$ is a SIN-group if and only if it is an extension
\beq\label{SIN-extension}
1\to N = V\times K\to G\to D\to 1,
\eeq
where $D$ is discrete, $V\times K$ is a direct product, $V\simeq \R^n$, $K$ is compact, and both $K$, $V$ are normal in $G$. We will use these notations further.

In section \ref{section_SIN}, we prove a duality theorem for Moore groups. These are groups whose all irreducible continuous unitary representations are finite-dimensional. Since the class of Moore groups is closed under passing to the quotients \cite[Vol.~II, p.~1452]{palmer}, $D$ in \eqref{SIN-extension} is a Moore group too. It is known that every discrete Moore group is a finite extension of an abelian group \cite[Theorem 12.4.26 and p.~1397]{palmer}.

The decomposition \eqref{SIN-extension} allows to use the Mackey's construction of induced representations \cite{mackey}. Mackey's machinery requires the group to be separable, but when inducing from an open subgroup, this is no longer needed.
Let $\pi$ be a unitary representation of $N$ in a Hilbert space $H$. There are then two realizations of the induced representation $T$. The first one acts on the space $L_2(G,H)$ of square-summable $H$-valued functions such that $f(\xi g) = \pi(\xi) f(g)$ for all $\xi\in N$, $g\in G$. The action is given by $T_g f(h) = f(hg)$.

In the second realization, $T$ acts on $L_2(D,H)$. We fix any map $s: D\to G$ such that $Ns(x) = x$, then $T$ acts as
\beq\label{induced_rep_formula}
(T_g f)(x) = \pi(s_x g s_{xg}^{-1}) f(xg).
\eeq

\blm\label{induced_rep}
If a unitary representation $\pi$ of\/ $N$ is norm continuous, then the induced representation $T$ of\/ $G$ is also norm continuous.
\elm
\bpr
For proof, we use the realization \eqref{induced_rep_formula} of $T$ on $L_2(D,H)$. Take $\e>0$. By assumption there is a neighborhood of identity $U$ in $N$, such that $\|\pi(g)-1_H\|<\e$ when $g\in U$. Since $N$ is open, $U$ is at the same time a neighborhood of identity in $G$. Since $G$ is a SIN-group, there exists another neighborhood of identity $W$ such that $gWg^{-1}\subset W\subset U$ for all $g\in G$. Take now $f\in L_2(D,H)$:
$$
\|T_gf-f\|^2 = \sum_{x\in D} \|\pi(s_x g s_{xg}^{-1}) (f(xg))-f(x)\|^2.
$$
If $g\in W\subset N$, then $xg=x$, so that
$$
\|T_gf-f\|^2 = \sum_{x\in D} \|\pi(s_x g s_x^{-1}) (f(x))-f(x)\|^2
 \le \sum_{x\in D} \|\pi(s_x g s_x^{-1})-1_H\|^2 \|f(x)\|^2.
$$
Moreover, $s_x g s_x^{-1}\in W$, so that $\|\pi(s_x g s_x^{-1})-1_H\|<\e$, and it follows that $\|T_gf-f\|<\e\|f\|$, what proves the lemma.
\epr

It is known that $N$ is a Moore group \cite[Theorems 12.4.16 and 12.4.28]{palmer}.

\blm\label{K(N)-tensor}
Linear combinations of the functions of the type $f(v)g(k)$, where $f\in\K(V)$, $g\in\K(K)$,
span a dense subspace in $\K(N)$.
\elm
\bpr
Let us show first that all functions of this type belong to $\mathcal K(N)$. By proposition \ref{K(G)-properties}(1), every function $f\in\mathcal K(V)$ is a coefficient of a norm continuous representation $\rho$ of $V$: $f(v)=\langle \rho(v)\xi_1,\eta_1\rangle$. We can define a norm continuous representation $\pi$ of $N$ as $\pi(vk)=\rho(v)$. Then $F(t)=\langle \pi(t)\xi_1,\eta_1\rangle$ will be a coefficient of $\pi$ such that $F(vk)=f(v)$, so $F\in\mathcal K(N)$.

Similarly, all functions of the type $G(vk)=g(k)$, $g\in\mathcal K(K)$, belong to $\mathcal K(N)$. Now by lemma \ref{K(G)-algebra}, $FG\in\K(N)$ as well, and $(FG)(vk) = f(v)g(k)$.

Next we need to prove that the linear span $E$ of such products is dense in $\K(N)$; due to lemma \ref{K(G)-properties}(5), it is sufficient to show that $E$ contains coefficients of every irreducible representation $\pi$ of $N$. Since $N$ is a Moore group, $\pi$ is finite-dimensional. Let $\pi$ act on a space $H$, and let $F(g)=\langle\pi(g)\xi,\eta\rangle$ be a coefficient of $\pi$. The restrictions $\rho=\pi|_V$, $\sigma=\pi|_K$ are representations of $V$ and $K$ respectively. Then $F$ may be represented in the following form, for $g=vk$, $v\in V$, $k\in K$:
$$
F(g)=F(vk)=\langle\pi(vk)\xi,\eta\rangle
 =\langle\rho(v)\sigma(k)\xi,\eta\rangle
 =\langle\sigma(k)\xi,\rho(v)^*\eta\rangle.
$$
Choose an orthonormal basis $\{e_i\}_{i=1}^n$ in $H$, then
$$
F(g) =\sum_{i=1}^n \langle\sigma(k)\xi,e_i\rangle \langle e_i,\rho(v)^*\eta\rangle
 = \sum_{i=1}^n f_i(v)g_i(k),
$$
where $f_i(v)=\langle \rho(v)e_i,\eta\rangle$ and $g_i(k)=\langle\sigma(k)\xi,e_i\rangle$ are coefficients of $\rho$ and $\sigma$ respectively. Thus, $F\in E$, and this proves the lemma.
\epr

\begin{nots}
By $I_X$, we denote the indicator function of a set $X$. By $\trans tf$ we denote the function $\trans tf(x)=f(xt^{-1})$.
\end{nots}

\blm\label{trans_K(N)}
Consider $\mathcal K(N)$ as a subset of $C(G)$, every function being extended by zero to $G\setminus N$. Then $I_{tN} \cdot\mathcal K(G) = \trans t{\mathcal K(N)}\equiv \{\trans tf: f\in\mathcal K(N)\}$ for every $t\in G$.
\elm
\bpr
This is enough to prove for $t\in N$, since for any $t\in G$ we would have
$$
\trans t{\mathcal K(N)} = \trans t{\big(I_N\cdot\mathcal K(G)\big)} =
\trans t{I_N} \cdot\trans t{\mathcal K(G)} = I_{Nt} \cdot\mathcal K(G) = I_{tN} \cdot\mathcal K(G).
$$

So let us prove first that $\mathcal K(N)\subset \mathcal K(G)$, this will imply that $\mathcal K(N) = I_N\cdot\mathcal K(N)\subset I_N\cdot\mathcal K(G)$. For every representation $\pi$ of $N$, acting in a space $H$, take the induced representation $T$ of $G$ (norm continuous if $\pi$ was so) acting on $L_2(D,H)$, see formula \eqref{induced_rep_formula}. Then for every function $f\in\mathcal K(N)$ of the type $f(g)=
\langle \pi(g)\xi, \eta\rangle$ with some $\xi,\eta\in H$ we have its zero-extension $F(g)=\langle T(g)(\xi I_N), \eta I_N\rangle$ which is in $\mathcal K(G)$.

Converse is even simpler: if $T$ is a representation of $G$, then its restriction onto $N$ is a representation of $N$, and so the restriction of every $f\in \mathcal K(G)$ to $N$ is a coefficient of a (norm-continuous) representation of $N$, which is in $\mathcal K(N)$. Since $I_N\cdot \mathcal K(G) = \{f|_N:f\in\mathcal K(G)\}$, we have proved the second inclusion and with it, the lemma.
\epr

\bcor\label{sep-points}
On any SIN-group $G$,
\begin{enumerate}
\item norm continuous unitary representations separate points;
\item the algebra $\mathcal K(G)$ separates points.
\end{enumerate}
\ecor

\bpr
Due to proposition \ref{K(G)-properties}, these statements are equivalent.
Since the subgroup $N$ is a Moore group, so finite-dimensional unitary representations separate its points. All of them are norm continuous. Now the statements follow directly from lemma \ref{trans_K(N)}.
\epr

\section{Structure of the dual algebra}\label{section-stereo}

Though we will use further only the case of Moore groups, the results of this section are naturally stated for all SIN-groups.

Theorems \ref{C(G)-abel} and \ref{direct_product_irr} give the description of the dual algebras for Abelian and compact groups. Now we describe the algebra $\hatC(G)$ for a discrete group $G$.

\btm\label{discrete}
A group $G$ is discrete if and only if\/ $\hatC(G)$ is a Banach algebra. In this case, $\hatC(G)$ is equal to the classical group \cst-algebra $C^*(G)$.
\etm
\bpr
Let $G$ be discrete. Then any representation of $G$ is norm continuous, thus all the \cst-seminorms are continuous on $M_c(G)$. The supremum $p_{\max}$ of all these seminorms is finite and equal to the classical \cst-norm of $\ell_1(G)$, restricted to $M_c(G)$ (which is in this case just the space of all finitely supported functions, so it is contained in $\ell_1(G)$). Clearly $M_c(G)$ is dense in $\ell_1(G)$, so its completion $C_{p_{\max}}^*(G)$ with respect to $p_{\max}$ is equal to that of $\ell_1(G)$, and the latter is the classical group \cst-algebra $C^*(G)$. From the other side, being a \cst-norm, $p_{\max}$ is contained in $\mathcal P(G)$. Since $p_{\max}$ is then the maximal element in $\mathcal P(G)$, the inverse limit $\hatC(G)=\liminv_{p\in\mathcal P(G)} C^*_p(G)$ equals to the corresponding Banach algebra $C_{p_{\max}}^*(G)$ \cite[Corollary 2.5.12]{engel}, so we come to the conclusion that $\hatC(G) = C^*(G)$.

Conversely, let $\hatC(G)$ be a Banach algebra. Then its \cst-norm $p\,$ belongs to $\mathcal P(M_c(G))$ and $\hatC(G)=C_p^*(G)$. We know that $G$ is homeomorphically imbedded into $M_c(G)$ via delta-functions (Proposition \ref{G-imbedded}), so it is continuously mapped to $\hatC(G)$. Thus, to prove that $G$ is discrete, it is sufficient to show that $p\,(\delta_t-\delta_e)\ge 1$ for any $t\ne e$. It is clear that $p\,(\delta_t-\delta_e)\ge |\phi(\delta_t)-\phi(\delta_e)|=|\phi(\delta_t)-1|$ for any state $\phi$ of any algebra $C_q^*(G)$, $q\in \mathcal P(G)$. If we consider this $\phi\in \big(C_q^*(G)\big)^\star\subset\mathcal K(G)$ as a function on $G$, it will be positive definite, and conversely, every positive-definite function $\phi\in\mathcal K(G)$ such that $\phi(e)=1$ is a state of some algebra $C_q^*(G)$. Since $\K(G)$ separates points of $G$ and is contained in the linear span of positive-definite functions (proposition \ref{K(G)-properties}), for the given $t\ne e$ we can always find a positive-definite function $\psi\in\K(G)$ such that $\psi(t)\ne \psi(e)=1$. We will have automatically $|\psi(t)|\le 1$ \cite[13.4.3]{dixmier}. By lemma \ref{K(G)-algebra}, $\K(G)$ is closed under pointwise multiplication and conjugation. If $|\psi(t)|=1$, then for some $n$ we have $|\psi(t)^n-1|\ge1$; if $|\psi(t)|<1$, then $\psi(t)^n\to0$, $n\to\infty$, so $|\psi(t)^n-1|\to1$. In both cases $\sup\limits_n |\psi(t)^n-1|\ge1$. Since $\psi^n$ is also a positive definite function and $\phi^n(e)=1$, this implies the estimate $p\,(\delta_t-\delta_e)\ge 1$.
\epr

Next we turn to direct products of SIN-groups and show (Corollary \ref{direct_product}) that they correspond to tensor products of group algebras. As with \cst-algebras (see, e.g., \cite{kad}), there is a range of possible tensor products on pro-\cst-algebras. We use the generalization of the maximal \cst-tensor product, which was defined for pro-\cst-algebras by N.~C.~Phillips \cite{phil}.  If $A=\liminv A_p$, $B=\liminv B_q$, then
\beq\label{AomaxB}
A\omax B = \liminv A_p\omax B_q
\eeq
over pairs $(p,q)$ directed coordinatewise \cite[Proposition 3.2]{phil}. This property, used sometimes without reference, may be taken as definition of the maximal tensor product $\omax$. As expected, $\omax$ is associative and commutative; $A\omax \C \simeq A$ for any $A$.

\btm\label{envelope_of_tensor_product}
Let $C$ be an algebra which is a locally convex space, and let $A$, $B$ be its commuting subalgebras such that the subalgebra $\<AB\>$ generated by $A$ and $B$ is dense in $C$. If every pair of continuous commuting representations of $A$ and $B$ in the same space may be extended to a continuous representation of\, $C$, then $C^\Diamond=A^\Diamond\omax B^\Diamond$.
\etm
\bpr
Let $A^\Diamond=\liminv A_p$, $B^\Diamond=\liminv B_q$, then $A^\Diamond\omax B^\Diamond= \liminv A_p\omax B_q$.
By a known property of the maximal \cst-tensor product \cite[11.3.4]{kad}, for every pair $(p,q)$ there are commuting representations $S$, $T$ of $A_p$ and $B_q$ such that $\|\gamma\|_{A_p\omax B_q}=\|(S\otimes T)(\gamma)\|$ for any $\gamma\in A_p\omax B_q$. We can consider $T$ and $S$ as representations of $A$ and $B$ respectively. The closure $Z_{pq}$ of the linear span of $S(A)T(B)$ is isomorphic to $A_p\omax B_q$. By assumption there is a representation $U$ of $C$ such that $U|_A=S$, $U|_B=T$. It is clear that the closure of $U(C)$ is also $Z_{pq}$. Let $r_{pq}$ be the corresponding \cst-seminorm on $C$, then the \cst-algebra $C_{r_{pq}}^*(C)$ is isomorphic to $Z_{pq}$. Thus, there exists an isomorphism $\phi_{pq}: A_p\omax B_q \to C_{r_{pq}}^*(C)$, which maps (see notations \ref{nots-envelope}) $i_p(a)\otimes i_q(b)$ to $i_{r_{pq}}(ab)$.

If we consider another pair of seminorms $p'\le p$, $q'\le q$ and define an isomorphism $\phi_{p'q'}$ in the same way, then the following diagram commutes:
$$
 \xymatrix @R=1.pc @C=2.pc
 {
 A_p\omax B_q \ar@{->}[rrr]^{\phi_{pq}}\ar@{->}[dd]^{}
 &
 & &
 C_{r_{pq}}^*(G) \ar@{->}[dd]^{}
 \\
 & & &
 \\
 A_{p'}\omax B_{q'}\ar@{->}[rrr]_{\phi_{p'q'}}
 & & &
 C_{r_{p'q'}}^*(G)
 }
$$

Thus, the maps $\phi_{pq}$ define a continuous map of inverse limits $\phi: A^\Diamond\omax B^\Diamond\to C^\Diamond$, provided that the set $\{r_{pq}: p\in\mathcal P(A), q\in\mathcal P(B)\}$ is cofinal in $\mathcal P(C)$ \cite[\S2.5]{engel}. Let us check the latter requirement.

For any $r\in\mathcal P(C)$ we have restrictions $p=r|_A\in\mathcal P(A)$ and $q=r|_B\in\mathcal P(B)$. By definition, the norm of every $\gamma\in A_p\omax B_q$ is equal to the supremum of $\|(S\otimes T)(\gamma)\|$ over all commuting representations $S$, $T$ of $A_p$ and $B_q$ respectively. We can consider $S=i_r|_A$ and $T=i_r|_B$ as commuting representations of $A$ and $B$ (since $A$ and $B$ commute in $C$). But $S$ and $T$ may be extended continuously to $A_p$ and $B_q$, since $\|S(a)\|=r(a)=p(a)$ and $\|T(b)\|=r(b)=q(b)$ for all $a\in A$, $b\in B$. Thus, $\|(S\otimes T)(\gamma)\|\le\|\gamma\|_{A_p\omax B_q}$ for any $\gamma\in A_p\omax B_q$. Now take $\gamma=(i_p\otimes i_q)(u)$ for some $u\in \langle AB\rangle$, then $\|(S\otimes T)(\gamma)\|=r(u)$ and $\|\gamma\|_{A_p\omax B_q}=r_{pq}(u)$, where $r_{pq}$ is defined as above. This, together with the density of $\langle AB\rangle$ in $C$, yields $r\le r_{pq}$. This means that $\mathcal P_{\max} = \{r_{pq}\}$ is a cofinal set in $\mathcal P(C)$.

So, we have a continuous map of inverse limits such that every coordinate map $\phi_{pq}$ is a homeomorphism. It follows \cite[Proposition 2.5.10]{engel} that this map is a homeomorphism of $A^\Diamond\omax B^\Diamond$ onto $C^\Diamond$. Obviously this is also a *-algebra homomorphism, so $A^\Diamond\omax B^\Diamond$ is isomorphic to $C^\Diamond$.
\epr

\bcor\label{direct_product}
$\hatC(G\times H)=\hatC(G)\omax \hatC(H)$ for any SIN-groups $G$ and $H$.
\ecor
\bpr
We can apply theorem \ref{envelope_of_tensor_product} with $A=M_c(G)$, $B=M_c(H)$ and $C=M_c(G\times H)$. Here $\<AB\>$ is the linear span of all delta-functions, since $\delta_{s,e}\cdot\delta_{e,t}=\delta_{s,t}$, and so it is dense in $M_c(G\times H)$. Next, every pair of commuting representations of $G$ and $H$ defines a representation of $G\times H$ (with the same type of continuity); and bijection between representations of groups and their measure algebras gives us the required extension property.
\epr

After several lemmas, we describe algebra $\hatC(G)$ in more detail. The main tool is Theorem \ref{decomp_for_hatC(G)}, which states that every norm continuous representation of $G$ is decomposed into a finite sum of representations such that each of them, restricted onto $K$, is a multiple of a single irreducible representation $\pi\in\widehat K$. It follows, see Corollary \ref{hatC_direct_prod}, that $\hatC(G)$ is a direct product of Fr\'echet pro-\cst-algebras $C^*_\pi(G)$ over all irreducible representations $\pi\in\widehat K$.

The first preparatory lemma is a particular case of \cite[Proposition~3.4]{phil}:

\blm\label{C(H)hatC(G)}
For any pro-\cst-algebra $A$ and a locally compact space $M$, $C(M)\omax A$ is isomorphic
to the algebra $C(M,A)$ of continuous functions from $M$ to $A$, in the topology of uniform convergence on compact sets.
\elm

\blm\label{p_t-continuous}
Let $G$ be a SIN-group. In notations \eqref{SIN-extension}, for $t\in D$ denote $M_t=\{\mu\in M_c(G):\supp \mu\subset tN\}$. Then a seminorm $p$ on $M_c(G)$ is continuous if and only if its restriction to every $M_t$, $t\in D$, is continuous.
\elm
\bpr
Consider the open subgroup $N$ in $G$. Since $G$ is the disjoint union of open cosets $tN$, $t\in D$, we have $C(G)=\prod_{t\in D} C(tN)$. Then the conjugate space $M_c(G)=\big(C(G)\big)^\star$ is the usual direct sum of $\big(C(tN)\big)^\star=M_c(tN)$ \cite[4.20]{akb}. We can identify $M_c(tN)$ with $M_t$, so $M_c(G)=\oplus_{t\in D} M_t$ (every measure $\mu\in M_c(G)$ is a finite sum of measures $\mu_t\in M_t)$. The base of neighbourhoods of zero in $M_c(G)$ is then the set of absolutely convex sets $U$ such that $U\cap M_t$ is a neighbourhood of zero in $M_t$ \cite[Proposition~V.2.4]{robertson}.

To prove the non-obvious implication of the lemma, suppose that $p$ is a seminorm on $M_c(G)$ such that its restriction to every $M_t$ is continuous. Then the set $U=\{\mu\in M_c(G): p(\mu)\le1\}$ is absolutely convex, and by assumption $U\cap M_t$ is an (open) neighbourhood of zero in $M_t$. It follows that $U$ is open, so that $p$ is continuous.
\epr

By lemma \ref{hom-of-hatC(G)}, the inclusion $N\to G$ generates an imbedding $\hatC(N)\hookrightarrow\hatC(G)$. Its image is the closure of the linear span of all $\delta_t$, $t\in N$, or equivalently the closure of the subspace of all $\mu\in M_c(G)$ such that $\supp \mu\subset N$.

\blm\label{sup_p_t}
Let $G$ be a SIN-group with a decomposition \eqref{SIN-extension}. Then
\begin{enumerate}
\item for any seminorm $p_0\in \mathcal P(N)$, there exists a seminorm $p\in\mathcal P(G)$ such that $p|_{\hatC(N)}=p_0$;
\item the supremum of all such seminorms $p_{\max}=\sup\{ p\in \mathcal P(G): p|_{\hatC(N)}=p_0 \}$ is finite and contained in $\mathcal P(G)$;
\item if\, $q\in\mathcal P(G)$ and\/ $q|_{\hatC(N)}\le p_0$, then $q\le p_{\max}$.
\end{enumerate}
\elm
\bpr
(1) Let $p_0\in \mathcal P(N)$. It is the norm of some representation of $N$, which induces a representation $T$ of $G$ (lemma \ref{induced_rep}); then the norm $p$ of $T$ is a \cst-seminorm in $\mathcal P(G)$ and $p|_{\hatC(N)}=p_0$.

(2) For any $\mu\in M_c(G)$ and $t\in D$, let $\mu_t$ be the restriction of $\mu$ onto $C(tN)$. Since $\supp \mu$ is compact and $N$ is open, only finite number of these $\mu_t$ is nonzero. For every~$t$, $\supp \delta_{t^{-1}}*\mu_t \subset N$, i.e. $\delta_{t^{-1}}*\mu_t \in\hatC(N)$. Given $p_0\in\mathcal P(N)$, these considerations allow to define a seminorm $Q$ on $M_c(G)$ as follows:
\beq
Q(\mu) = \sum_{t\in D} p_0( \delta_{t^{-1}}*\mu_t).
\eeq
First of all, this sum is finite for any $\mu\in M_c(G)$, so it is well defined. Next, the restriction of $Q$ onto every $M_t$ (in notations of lemma \ref{p_t-continuous}) is the composition of $p_0$ and a translation by $\delta_{t^{-1}}$, so it is continuous. By lemma \ref{p_t-continuous}, $Q$ is continuous on $M_c(G)$. Of course, it need not in general be a \cst-seminorm.

If now $p$ is any \cst-seminorm on $M_c(G)$, then $p(u*\mu)=p(\mu*u)=p(\mu)$ for any unitary element $u\in M_c(G)$ and any $\mu\in M_c(G)$. 
In particular, $p(\delta_t*\mu)=p(\mu)$ for all $t\in G$, $\mu\in M_c(G)$.

Let now $p|_{\hatC(N)}=p_0$. Using notations introduced above, we have now:
\beq
p(\mu) \le \sum_{t\in D} p(\mu_t) = \sum_{t\in D} p( \delta_{t^{-1}}*\mu_t) = \sum_{t\in D} p_0( \delta_{t^{-1}}*\mu_t)  = Q(\mu).
\eeq
Thus, $p\le Q$, so $p_{\max}=\sup\{p\in\mathcal P(G): p_{\hatC(N)}=p_0\}\le Q<\infty$. Since $p_{\max}$ is bounded by a continuous seminorm on $M_c(G)$, it is itself continuous. And obviously it has the \cst-property since it is the supremum of a family of \cst-seminorms.

(3) Let $r=\max(q,p_{\max})$. This is a \cst-seminorm, and $r|_{\hatC(N)}= p_0$; so, by (2) $r\le p_{\max}$. By definition of $p_{\max}$, we have $r=p_{\max}$, that is, $q\le p_{\max}$.
\epr

\blm
Let $\pi\in\widehat K$ and let $\chi_\pi$ be the character of $\pi$. Define a measure $\tilde\chi_\pi$ on $G$ by
$\tilde\chi_\pi(f)=\int_K f\chi_\pi$, $f\in C(G)$. Then $\tilde\chi_\pi$ is central in $M_c(G)$.
\elm
\bpr
Since the linear span of delta-functions is dense in $M_c(G)$, it is sufficient to prove that $\tilde\chi_\pi$ commutes with $\delta_t$ for all $t\in G$. For $t\in G$, let $\tilde\xi_t=\delta_t*\tilde\chi_\pi*\delta_t^{-1}$; we need to show that $\tilde\xi_t=\tilde\chi_\pi$. By direct calculation, using the normality of $K$, we get that $\tilde\xi_t(f)=\int_K f(u)\chi_\pi(t^{-1}ut)du$ for all $f\in C(G)$. Denote $\xi_t(u)=\chi_\pi(t^{-1}ut)$, then the identity $\tilde \xi_t=\tilde\chi_\pi$ is equivalent to $\xi_t=\chi_\pi$, which we will continue to prove.

First, $\xi_t$ is central in $M_c(K)$, i.e. $\xi_t(s^{-1}us)=\xi_t(u)$ for all $s,u\in K$. To prove this, take any $s\in K$ and choose (by normality of $K$) such $s'\in K$ that $s't=ts$. We know that $\chi_\pi$ is central as every character; then $\xi_t(s^{-1}us) = \chi_\pi(t^{-1}s^{-1}ust)=\chi_\pi(s'^{-1}t^{-1}uts')=\chi_\pi(t^{-1}ut)=\xi_t(u)$.

It is known that then $\xi_t=\sum_{\sigma\in\widehat K} \lambda_{t,\sigma}\chi_\sigma$ with some $\lambda_{t,\sigma}\in\C$, the series converging in $L_2(K)$, and
$$
\sigma(\tilde\xi_t)=\int_K \xi_t(x)\sigma(x)dx= (\dim\sigma)^{-1}\lambda_{t,\sigma}\Id
$$
for every $\sigma\in\widehat K$ \cite[28.49, 28.50]{HR}. If $\sigma\ne\pi$, then (see lemma \ref{sup_p_t})
$$
\|\sigma(\tilde\xi_t)\|=p_\sigma^{\max}(\tilde\xi_t)=p_\sigma^{\max}(\delta_t*\chi_\pi*\delta_t^{-1})= p_\sigma^{\max}(\chi_\pi)=\|\sigma(\chi_\pi)\|=0,
$$
so $\lambda_{t,\sigma}=0$. Thus, $\xi_t=\lambda_{t,\pi}\chi_\pi$. Finally, from $\xi_t(e) = \chi_\pi(t^{-1}et)=\chi_\pi(e)$ we deduce that $\lambda_{t,\pi}=1$ and $\xi_t=\chi_\pi$, what proves the lemma.
\epr

\btm\label{decomp_for_hatC(G)}
Every norm continuous representation $T$ of $G$ is a finite direct sum of representations $T_\pi$ such that $T_\pi|_K$ is a multiple of a single irreducible representation $\pi\in\widehat K$.
\etm
\bpr
For $\pi\in \widehat K$, let $d_\pi=\dim\pi$. In notations of the previous lemma, define $T_\pi(\mu) = d_\pi T(\tilde\chi_\pi \mu)$ for all $\mu\in M_c(G)$. It is known \cite[27.24]{HR} that $\tilde\chi_\pi \tilde\chi_\sigma=d_\pi^{-1}\tilde\chi_\pi$ if $\pi=\sigma$ and 0 otherwise. This implies, together with centrality of $\tilde\chi_\pi$, that
$$
T_\pi(\mu\nu) = d_\pi T(\tilde\chi_\pi \mu \nu)
= d_\pi T(d_\pi\tilde\chi_\pi^2 \,\mu  \nu )
= d_\pi^2 T(\tilde\chi_\pi \mu ) T( \tilde\chi_\pi \nu)
=T_\pi(\mu)T_\pi(\nu),
$$
so $T_\pi$ is a representation of $G$. Clearly it is a subrepresentation of $T$ on the image of $\tilde\chi_\pi$. If we consider only restrictions to $K$, then $T_\pi$ is known to be a multiple of $\pi$, and the images of $T_\pi$ are orthogonal for different $\pi\in\widehat K$ \cite[27.44]{HR}. Thus, $T=\oplus_{\pi\in\widehat K} T_\pi$, and from lemma \ref{compact-finite} it follows that $T_\pi$ is nonzero only for a finite number of $\pi$.
\epr

Before we can formulate corollaries on the structure of $\hatC(G)$, let us describe the seminorms on $\hatC(N)$. In notations~\eqref{SIN-extension}, by lemmas \ref{direct_product}, \ref{C(G)-abel} and \ref{C(H)hatC(G)}, we have
\beq\label{N_e}
\hatC(N)= \hatC(\R^n\times K) = \hatC(\R^n)\!\omax\! \hatC(K)
 = C(\R^n)\! \omax\! \hatC(K) = C(\R^n, \hatC(K)).
\eeq
Thus, every $\mu\in\hatC(N)$ may be considered as a continuous function from $\R^n$ to $\hatC(K)=\prod_{\pi\in\widehat K}C^*_\pi(K)$.
A defining system of seminorms on $\hatC(N)$
has then the following form: for $k\in\N$ and $\sigma=\{\pi_1,\ldots,\pi_m\}\subset\widehat K$
\beq\label{p_on_N}
p_{k,\sigma}(\mu) = \sup_{|x|\le k} \;\max_{\pi\in\sigma} \|\mu(x)\|_{C^*_\pi(K)},
\eeq
where $|x|$ is any fixed norm in $\R^n$.

\bprop\label{order_for_p_k,sigma}
For every $\sigma\subset \widehat K$ and $k\in\N$, $p_{k,\sigma}^{\max} = \max\{ p_{k,\pi}^{\max}: \pi\in\sigma\}$.
\eprop
\bpr
It is clear that on $\hatC(N)$, we have $p_{k,\pi}\le p_{k,\sigma}$ for all $\pi\in\sigma$. By lemma \ref{sup_p_t}(3), $p_{k,\pi}^{\max}\le p_{k,\sigma}^{\max}$ and $p_{k,\sigma}^{\max} \ge \max\{ p_{k,\pi}^{\max}: \pi\in\sigma\}$.

To prove the converse inequality, choose a representation $T$ such that $p_{k,\sigma}^{\max}=\|T\|$. By the preceding theorem we have direct sum decomposition $T=\oplus_{\pi\in s} T_\pi$ over some finite set $s\subset \widehat K$.

For every $\pi\in s$, $T_\pi|_N$ is a direct product of a representation $\rho$ of $V\simeq\R^n$ and of $T_\pi|_K$ which is a multiple of $\pi$. Since we are interested only in the norm of $T_\pi|_N$, we can assume that $T_\pi|_K=\pi$. For $\rho$, there is a compact set $M\subset\R^n$ so that $\rho(t)(x)=e^{itx}$ for $t\in V$, $x\in M$. Then $T_\pi|_N$ can be written as $T(t,\kappa)(x) =e^{itx}\pi(\kappa)$ for $t\in V$, $\kappa\in K$, $x\in M$.

Since $\|T_\pi|_N\|\le\|T|_N\|=p_{k,\sigma}$, it follows that $\rho(t)(x)=0$ if $|x|>k$, so $M\subset \{x:|x|\le k\}$. But then $\|T_\pi|_N\|\le p_{k,\pi}$.

Now by lemma \ref{sup_p_t}(3) $\|T_\pi\| \le p^{\max}_{k,\pi}$, and then
$$
p^{\max}_{k,\sigma}=\|T\|=\|\oplus_{\pi\in s}T_\pi\| = \max_{\pi\in s} \|T_\pi\| \le \max_{\pi\in s} p^{\max}_{k,\pi},
$$
what proves the statement.
\epr

\bcor\label{hatC_direct_prod}
For a SIN-group $G$, $\hatC(G)=\prod_{\pi\in\widehat K} C_\pi$ where $C_\pi=\liminv_k C^*_{p_{k,\pi}^{\max}}(G)$ if $V\ne\{0\}$ and $C_\pi=C^*_{p_\pi^{\max}}(G)$ otherwise.
\ecor
\bpr
First of all, the set $\mathcal P_{\max}(G)=\{p_{k,\sigma}^{\max}: k\in\N, \sigma\subset \widehat K\}$ is cofinal in $\mathcal P(G)$ by lemma \ref{sup_p_t}. Proposition \ref{order_for_p_k,sigma} shows that the order on $\mathcal P_{\max}(G)$ is the natural order on $\N$ and inclusion of subsets of $\widehat K$.  Thus,
$$
\hatC(G) = \liminv_{k\in\N,\sigma\subset\widehat K} \prod_{\pi\in\sigma} C^*_{p^{\max}_{k,\pi}}(G).
$$
From the other side, using \eqref{AomaxB}, we have
\begin{align*}
\prod_{\pi\in\widehat K} C_\pi &= \liminv_{\sigma\subset \widehat K} \prod_{\pi\in\sigma} C_\pi
= \liminv_{\sigma\subset \widehat K} \prod_{\pi\in\sigma} \liminv_{k\in\N} C^*_{p^{\max}_{k,\pi}}(G)
= \liminv_{\sigma\subset \widehat K} \Big(\liminv_{k\in\N} \prod_{\pi\in\sigma} C^*_{p^{\max}_{k,\pi}}(G)\Big).
\end{align*}
Now it is easy to see that the natural map between these two inverse limits is an isomorphism.
\epr

\bcor\label{PF}
The category $\mathfrak{PF}$ of locally convex spaces which are direct products of Fr\'echet spaces, contains $C(G)$ and $\hatC(G)$ for all SIN-groups $G$. Every space in this category is stereotype. If $A,B\in\mathfrak{PF}$ are direct products of Fr\'echet pro-\cst-algebras, then $A\omax B\in\mathfrak{PF}$.
\ecor
\bpr
Due to decomposition \eqref{SIN-extension}, $C(G)$ is homeomorphic to $C(N)^D$ as a locally convex space. Every $C(N)$ is a Fr\'echet space because $N$ is sigma-compact; thus, $C(G)\in\mathfrak{PF}$. Every $C_\pi$ in lemma \ref{hatC_direct_prod} is a countable inverse limit of Banach spaces, so it is a Fr\'echet space, thus $\hatC(G)\in\mathfrak{PF}$. Direct product of Fr\'echet spaces is stereotype by \cite[4.4 and 4.20]{akb}.

Next, if $A=\prod A_\a$ and $B=\prod B_\beta$ where every $A_\a$, $B_\beta$ is a Fr\'echet pro-\cst-algebra, then $A\omax B=\prod (A_\a\omax B_\beta)$, and every $A_\a\omax B_\beta$ is Fr\'echet. This shows that $\mathfrak{PF}$ is closed under maximal tensor products.
\epr

\section{Strict algebras and Moore groups}\label{section-Moore}

In this section we consider a category $\mathfrak{Moore}$ which includes the algebras $C(G)$, $\hatC(G)$ for every Moore group $G$ (a group is called {\it Moore} if all its irreducible representations are finite-dimensional). This category is consists of Hopf algebras with respect to a topological tensor product $\odot$, satisfying the algebraic definition of a Hopf algebra in a tensor category.

Recall the following general definitions of an algebra, coalgebra and a Hopf algebra in a tensor category. We will write $\C$ for the unit object, because in our categories it is the space of the complex numbers.

\begin{definition}\label{def-hopf}
An {\it algebra} in a tensor category with the tensor product $\otimes$ and the unit object $\C$ is an object $A$ with morphisms $\m:A\otimes A\to A$ (multiplication) and $\imath: \C\to A$ (unit), such that\\
\null\qquad(alg.1) \ $\m(\m\otimes\id)=\m(\id\otimes\m)$,\\
\null\qquad(alg.2) \ $\m(\id\otimes\imath)=\id$, \ $\m(\imath\otimes\id)=\id$.

A {\it coalgebra} is an object $A$ with morphisms $\Delta:A\to A\otimes A$ (comultiplication) and $\e:A\to\C$ (counit) such that\\
\null\qquad(coal.1) \ $(\Delta\otimes\id)\Delta=(\id\otimes\Delta)\Delta$,\\
\null\qquad(coal.2) \ $(\id\otimes\e)\Delta=\id$, \ $(\e\otimes\id)\Delta=\id$.

Recall that for every algebra $A$, $A\otimes A$ is also an algebra. A {\it Hopf algebra} is an object $A$ which is an algebra and a coalgebra with an additional morphism $S:A\to A$ (antipode) such that\\
\null\qquad(hopf.1) $\Delta$ and $\e$ are algebra homomorphisms, i.e.\\
\null\hskip1.7cm $\Delta\m = \m_{A\otimes A}(\Delta\otimes \Delta)$ and $\e\,\m=\m_{\C\otimes \C}(\e\otimes\e)$.\\
\null\qquad(hopf.2) $\m(S\otimes\id)\Delta=\m(\id\otimes S)\Delta=\e\imath$.
\end{definition}

To indicate the tensor product we use the terms $\otimes$-algebra, $\otimes$-coalgebra and $\otimes$-Hopf algebra. In this section we consider the category $\mathfrak{Ste}$ of stereotype locally convex spaces with two tensor products: injective $\odot$ and projective $\circledast$. These tensor products are related by the identity $(X\circledast Y)^\star = X^\star\odot Y^\star$ \cite[Theorem~7.6]{akb}. $\mathfrak{Ste}$ is a tensor category with both \cite[Theorem~7.7]{akb}. As expected, if $A$ is a $\odot$-algebra ($\circledast$-algebra), then $A\odot A$ is also a $\odot$-algebra (respectively $A\circledast A$ is a $\circledast$-algebra) \cite[Theorems 10.15, 10.22]{akb}. For any locally compact group $G$, $C(G)$ is a $\odot$-Hopf *-algebra, and $M_c(G)$ is a $\circledast$-Hopf *-algebra \cite[10.24]{akb}.

If $A$ is a linear topological space and a $\otimes$-algebra, in analytical language this means that its multiplication $\m: A\times A\to A$ may be extended to a continuous map $\m: A\otimes A\to A$.
The notations $\m$, $\imath$, $\Delta$, $S$ and $\e$ will be used by default for the multiplication, unit, comultiplication, antipode and counit in any algebra or coalgebra.

\begin{definition}\label{def-star-hopf}
An {\it involution} on an algebra $A$ is an anti-homomorphism $*:A\to A$ such that $**=\id$. If $A$ is an algebra with involution then so is $A\otimes A$, with involution $*\otimes*$. An algebra morphism $\phi$ is a *-homomorphism if $\phi(a)^*=\phi(a^*)$. A {\it $\otimes$-Hopf *-algebra} is a $\otimes$-Hopf algebra with involution such that $\Delta$ and $\e$ are *-homomorphisms.
\end{definition}

It is clear that with the usual involution $C(G)$ is a $\odot$-Hopf *-algebra and $M_c(G)$ is a $\circledast$-Hopf *-algebra. This is a particular case of a general theorem \cite[10.23]{akb}: the conjugate space $A^\star$ to a $\odot$-Hopf algebra $A$ is a $\circledast$-Hopf *-algebra, and vice versa. Recall that the structural maps on $A^\star$ are given by conjugation: $\m_{A^\star}=\Delta_A^*$, $\Delta_{A^\star}=\m_A^*$, $\imath_{A^\star}=\e_A^*$, $\e_{A^\star}=\imath_A^*$, $S_{A^\star}=S_A^*$, and
\beq\label{Astar_involution}
\mu^*(a) = \overline{\mu\big( (Sa)^* \big)}
\eeq
for all $\mu\in A^\star$, $a\in A$.

In principle, the stereotype injective tensor product $\odot$ differs from the injective tensor product $\check\otimes$ in the usual theory of locally convex spaces \cite[\S~VII.2]{robertson}. But they coincide for Fr\'echet spaces, if one of the spaces has the approximation property (AP for short), and for their direct products \cite[Theorem~7.21 and the last formula in the proof of Lemma~8.6]{akb}. This class of spaces may be equivalently described as spaces with the AP which are direct products of Fr\'echet spaces (i.e. are in category $\mathfrak{PF}$, see corollary \ref{PF}), since AP is inherited by complemented subspaces \cite[p.~61]{cooper}. In particular, these tensor products coincide if both factors are $C(G)$ where $G$ is a SIN-group, since 1)~$C(G)\in \mathfrak{PF}$; 2)~AP is preserved by projective limits of Banach spaces \cite[p.~59]{cooper}, and $C(K)$ for a compact space $K$ has AP.

This allows to establish a link with the theory of strict Banach algebras. This notion was first introduced by N.~Varopoulos under the term {\it injective algebras}.

\begin{definition}
A Banach algebra $A$ is called {\it strict} if its multiplication can be extended to a continuous linear map from  $A\check\otimes A$ to $A$ (where $\check\otimes$ is the injective tensor product).
\end{definition}

For \cst-algebras, this property is equivalent to the following one \cite[Theorem 6]{aristov}:

\begin{definition}
A \cst-algebra $A$ is said to be of {\it bounded degree\/} if there exists a natural number $n$ such that all irreducible representations of $A$ are finite-dimensional and their dimensions do not exceed $n$.
\end{definition}

O.~Aristov had proved in \cite[Theorem 6]{aristov} that a \cst-algebra is strict if and only if it is of bounded degree, and in \cite[Theorem~7]{aristov} that this is also equivalent to topological isomorphism of $A\mathop{\check\otimes} A$ and $A\omin A$, where $\underset{\min}\otimes$ is the minimal \cst-tensor product. There are many other equivalent properties of \cst-algebras of bounded degree \cite[Theorem~2.5]{lau}. Now the following theorem follows:

\btm\label{strict-is-odot-algebra}
Let $A=\liminv A_p$ be a stereotype pro-\cst-algebra, $A\in\mathfrak{PF}$, and let every $A_p$ be a strict \cst-algebra. Then $A$ is a $\odot$-algebra.
\etm
\bpr
Since $\check\otimes$ commutes with inverse limits, $A\mathop{\check\otimes}A=\liminv A_p\mathop{\check\otimes}A_q$ over $p,q\in\mathcal P(A)$. Every pair $(p,q)$ is majorated by either $(p,p)$ or $(q,q)$, so the set $(p,p)$, $p\in\mathcal P(A)$, is cofinal and hence $A\mathop{\check\otimes}A=\liminv A_p\mathop{\check\otimes}A_p$. For every $p$, multiplication is a continuous map $\m_p:A_p\mathop{\check\otimes}A_p\to A_p$. These maps clearly agree with the inverse structure, so they define a continuous map of the inverse limits $\m: A\mathop{\check\otimes}A\to A$. Obviously $\m(a\otimes b)=ab$, so we have proved that $A$ is a $\mathop{\check\otimes}$-algebra.

Since every $A_p$ is of bounded degree, it is nuclear; then it has the approximation property \cite[p.~59]{cooper}. Since projective limits of Banach spaces preserve the approximation property \cite[p.~59]{cooper}, $A$ has this property too. It follows that $A\mathop{\check\otimes}A=A\odot A$, and the theorem is proved.
\epr

\begin{definition}
Let $A=\liminv A_p$ be a pro-\cst-algebra and let every $A_p$ be a strict \cst-algebra. Then $A$ is called a {\it strict pro-\cst-algebra}.
\end{definition}

These tensor products also coincide with the maximal tensor product:
\blm
If $A,B\in\mathfrak{PF}$ are pro-\cst-algebras and one of them is strict then
$A\odot B=A\omax B$.
\elm
\bpr
Similarly to the proof of theorem \ref{strict-is-odot-algebra}, we can show that $A\odot B=A\mathop{\check\otimes}B=\liminv A_p\mathop{\check\otimes}B_q$. By \cite[Theorem~7]{aristov}, $A_p\mathop{\check\otimes} B_q=A_p\omin B_q$. Since a strict \cst-algebra is nuclear, we have also $A_p\omin B_q = A_p\omax B_q$. Finally,
\begin{align*}
A\odot B = \liminv A_p\mathop{\check\otimes} A_q=\liminv A_p\omax A_q = A\omax A.
\end{align*}
\epr

\blm\label{degree}
If $G$ is a Moore group and $p\in {\cal P}(G)$, then the algebra $C_p^*(G)=C_p^*(M_c(G))$ (defined in \ref{nots-envelope}) is of bounded degree.
\elm
\bpr
I. It is known \cite[Theorem 12.4.27]{palmer} that a Moore group is equal to the projective limit of Lie Moore groups. By \cite[Theorem~1]{shtern}, every norm continuous representation of $G$ may be factored through one of the quotient groups in this limit, so we can assume that $G$ is itself a Lie group.

II. By \cite[Theorem 12.4.27]{palmer}, $G$ is a finite extension of a central group $G_1$ (a group is called central if its quotient over the center is compact). By \cite[Theorem 4.4]{gros-mosk-central} $G_1$ is equal to the direct product $\R^m \times G_2$, where $G_2$ has a compact open normal subgroup $K$. Let $Z$ be the center of $G_1$, and let $H=ZK$. Then $H$ is normal in $G_1$ and contains $Z$; it also contains $\R^m \times K$, so it is open in $G_1$. Thus $G_1/H$ is discrete and compact, hence finite. As $G/G_1 = (G/H) / (G_1/H)$ is also finite, we see that $G/H=F$ is finite too.

Suppose we have proved the theorem for the subgroup $H$. Take any $p\in {\cal P}(G)$ and consider the algebra $C_p^*(G)$. If $\tau$ is an irreducible representation of $C_p^*(G)$, then the corresponding representation of $G$ is also irreducible; since $G$ is a Moore group, $\tau$ is finite-dimensional. Denote $|F|=m$. By \cite[Theorem 1]{clifford}, the restriction of $\tau$ onto $H$ decomposes into at most $m$ irreducible representations of $H$.
Every component $\tau_i$ in this decomposition is continuous with respect to $p$ and therefore generates a representation of $C^*_p(H)$, which is also irreducible. By assumption $\dim\tau_i$ is at most $n=\deg C_p^*(H)$, so $\dim \tau\le mn$, what proves the lemma (assumed we have proved it for $H$).

III. It remains now to prove the theorem for the subgroup $H$, that is, a group representable as $H=ZK$, where the subgroup $Z$ is central and $K$ is compact and normal. Let $\pi: H\to {\cal B}(E)$ be a (norm-continuous) representation of $H$ on a space $E$. Let $\rho = \pi|_K$, $\sigma = \pi|_Z$, and let $C_\pi$, $C_\rho$, $C_\sigma$ be the \cst-algebras generated by $\pi$ and by these restrictions $\rho$, $\sigma$ respectively.

Since $\rho$ and $\sigma$ commute, $C_\pi$ is a continuous image of the maximal tensor product $C_\rho\omax C_\sigma$, so every irreducible representation of $C_\pi$ generates an irreducible representation of $C_\rho\omax C_\sigma$. This tensor product may be described explicitly. As it is shown in theorem \ref{C(G)-abel}, $C_\sigma$ is isomorphic to the algebra $C(M)$ for some compact set $M\subset \widehat Z$. By lemma \ref{compact-finite}, $C_\rho$ is isomorphic to a finite direct sum of matrix algebras: $C_\rho = \oplus_{j=1}^m M_{n_j}(\C)$. Thus,
$$
C_\rho\omax C_\sigma = \mathop{\oplus}_{j=1}^m C(M)\otimes M_{n_j}(\C)
$$
(here $\otimes$ is just the algebraic tensor product). For every $j$, a representation of $C(M)\otimes M_{n_j}(\C)$ is irreducible if and only if its restriction onto $1\otimes M_{n_j}(\C)$ is irreducible. Thus, the dimension of such a representation is at most $n_j^2$. Since the summands in this sum commute, it is easy to see that any irreducible representation of $C_\rho\omax C_\sigma$ has dimension at most $n_1^2\cdots n_m^2$, and this concludes the proof.
\epr

\bcor\label{moore-mult}
If $G$ is a Moore group, then $\hatC(G)$ is a $\odot$-algebra.
\ecor

\bprop\label{nomult}
For a discrete group $G$ which is not a Moore group, $\hatC(G)$ is not an $\odot$-algebra.
\eprop
\bpr
According to theorem \ref{discrete}, $\hatC(G)$ is the Banach algebra $C^*(G)$. From the proof of \cite[Theorem~7.21]{akb} one can see that $C^*(G)\mathop{\check\otimes}C^*(G)$ is imbedded into $C^*(G)\odot C^*(G)$ as a dense subspace. Thus, if $C^*(G)$  is an $\odot$-algebra, it is also a $\mathop{\check\otimes}$-algebra, i.e. it is strict. It follows that it is of bounded degree, so the dimensions of all irreducible representations of $G$ must be bounded by a common constant; a fortiori, they are finite-dimensional, so $G$ is a Moore group.
\epr

In the following lemma we prove that $A^{\star\Diamond}$ is automatically a $\odot$-Hopf *-algebra under certain conditions, which hold, in particular, for $A=C(G)$ on a Moore group $G$ (see further corollary~\ref{hatC-moore-hopf}).
One of the conditions is that the antipode $S$ is a *-antihomomorphism, i.e. $S(a^*)=S(a)^*$ (besides the usual condition $S(ab)=S(b)S(a)$). This is equivalent to the identity $S^2={\rm id}$ 
\cite[4.1]{chari}, and holds, in particular, for every commutative or cocommutative Hopf algebra.

\blm\label{moore-envelope}
Let $A$ be a $\odot$-Hopf *-algebra such that\\
1) its antipode is a *-antihomomorphism;\\
2) $\widehat A=A^{\star\Diamond}$ is a strict pro-\cst-algebra;\\
3) $(A^\star\circledast A^\star)^\Diamond = \widehat A\omax \widehat A$.\\
Then $\widehat A$ with structural maps extended by continuity from $A^\star$ is a $\odot$-Hopf *-algebra.
\elm
\bpr
By theorem \cite[10.23]{akb}, $A^\star$ is a $\circledast$-Hopf algebra. Immediate calculation shows that $\Delta_{A^\star}$, $\e_{A^\star}$ and $S_{A^\star}$ commute with involution. Thus, $\Delta_{A^\star}$ and $\e_{A^\star}$ are *-ho\-mo\-morph\-isms and by the universality property they are extended to *-ho\-mo\-morph\-isms of respective \cst-envelopes.

To show that $S_{A^\star}$ is continuous in the topology of $\widehat A$, take a seminorm $p\in\cal P(A^\star)$ and denote $p_S(a)=p(S_{A^\star}a)$. Then $p_S$ is a seminorm, and
$$
p_S(a^*a) = p(S_{A^\star}\!(a^*a)) = p(S_{A^\star}\!(a) S_{A^\star}\!(a^*)) = p(S_{A^\star}\!(a) S_{A^\star}\!(a)^*)
 = p(S_{A^\star}\!(a))^2 = p_S(a)^2,
$$
so $p_S\in \cal P(A^\star)$. This implies that $S_{A^\star}$ is continuous in the topology generated by all $p\in\mathcal P(A^\star)$, and this is the topology of $\widehat A$. So it can be extended by continuity to $\widehat A$.

The axioms (alg.1,2), (coal.1,2) and (hopf.1) hold obviously. For (hopf.2), note that $\widehat A\omax \widehat A = \widehat A\odot \widehat A$, and by assumption $\widehat A$ is a $\odot$-algebra. This implies that the morphisms $\m_{\widehat A}$, $\id\otimes S_{\widehat A}$, $S_{\widehat A}\otimes\id$ are well defined (on $\widehat A\odot \widehat A$), and then by density we get the identity (hopf.2).
\epr

\bcor\label{hatC-moore-hopf}
For any Moore group $G$, \;$\hatC(G)$ is a $\odot$-Hopf *-algebra.
\ecor
\bpr
We can apply lemma \ref{moore-envelope} with $A=C(G)$, because $C(G)\odot C(G)=C(G\times G)$ \cite[Theorem~8.4]{akb}, and with theorem \ref{direct_product} we have for $A^\star=M_c(G)$:
\begin{align*}
\big(M_c(G)\circledast M_c(G)\big)^\Diamond &= \big( C(G)\odot C(G) \big)^{\star\Diamond}
= \big( C(G\times G) \big)^{\star\Diamond}
= M_c(G\times G)^\Diamond \\
&= \hatC(G\times G) = \hatC(G)\omax \hatC(G) =
 M_c(G)^\Diamond\omax M_c(G)^\Diamond.
\end{align*}
\epr

\begin{definition}
Let $\mathfrak {Moore}$ denote the subcategory of Hopf $\odot$-algebras $A$ such that:
\\(i) the dual algebra $\widehat A$ is well defined;
\\(ii) $\hathatA=A$.
\\
Morphisms in $\mathfrak {Moore}$ are morphisms of $\odot$-Hopf *-algebras.
\end{definition}

\btm\label{functor_Moore}
$A\mapsto\widehat A$ is a contravariant functor on $\mathfrak {Moore}$.
\etm
\bpr
If $\phi:A\to B$ is a morphism in $\mathfrak {Moore}$, then $\widehat \phi$ is defined in the following natural way. The conjugate map $\phi^\star: B^\star\to A^\star$ is a morphism of $\circledast$-Hopf *-algebras, so in particular this is a *-homomorphism. By universality it may be extended to a *-homomorphism of respective \cst-envelopes $\widehat \phi: \widehat B\to \widehat A$.

By \cite[\S7.3, \S7.4]{akb} there are natural continuous maps $i_A: A^\star\circledast A^\star\to \widehat A\odot\widehat A$ and $i_B: B^\star\circledast B^\star\to \widehat B\odot\widehat B$, identical on elementary tensors. By assumption Hopf structures on $\widehat A$ and $\widehat B$ are well defined, what means that there are comultiplications $\Delta_{\widehat A}:\widehat A\to \widehat A\odot \widehat A$ and $\Delta_{\widehat B}: \widehat B\to \widehat B\odot\widehat B$. On $A^\star$ and $B^\star$ respectively, we have $\Delta_{\widehat A}=i_A \circ\Delta_{A^\star}$ and $\Delta_{\widehat B}=i_B\circ \Delta_{B^\star}$.

Consider the maps $\Delta_{\widehat A}\circ\widehat\phi$ and $(\widehat\phi\odot\widehat\phi)\circ\Delta_{\widehat B}$ from $\widehat B$ to $\widehat A\odot\widehat A$. On the dense subspace $B^\star\subset \widehat B$, they equal respectively $i_A\circ\Delta_{A^\star}\circ\phi^*$ and $i_A\circ(\phi^*\circledast\phi^*)\circ\Delta_{B^\star}$ (the latter holds since $\Delta_{B^\star}(B^\star)\subset B^\star\circledast B^\star$). But these maps are equal since $\phi^*$ is a morphism of $\circledast$-coalgebras. Thus, $\Delta_{\widehat A}\circ\widehat\phi=(\widehat\phi\odot\widehat\phi)\circ\Delta_{\widehat B}$, i.e. $\widehat \phi$ is a morphism of $\odot$-coalgebras.

If $\xi=\psi\circ\phi$, then $\xi^*=\phi^*\circ\psi^*$, so $\widehat\xi=\widehat\phi\circ\widehat\psi$. This proves that the functor $\widehat\;\;$ is contravariant. And finally, since $\hathatA=A$, this is a duality on $\mathfrak {Moore}$.
\epr

In section \ref{section_SIN} we will show that every algebra $C(G)$ and $\hatC(G)$, $G$ being a Moore group, is contained in $\mathfrak {Moore}$.

\section{Duality theorems}\label{section_SIN}

In this section we prove the main theorem \ref{SIN-reflexive}. The key task is to describe the space of characters of $\K(G)$, for a Moore group $G$. Then it follows that the dual algebra of $\hatC(G)$ is $C(G)$, so that $C(G)$ is a reflexive $\odot$-Hopf *-algebra.

\blm\label{K(G)-compact}
Let $G$ be a compact group. For every nonzero involutive character $\chi$ of $\K(G)$ there is $t\in G$ such that $\chi(f)=f(t)$ for all $f\in \K(G)$.
\elm
\bpr
Among the representations corresponding to semi\-norms $p$ there are, in particular, all irreducible unitary continuous representations. It follows that the algebra $\K(G)$ contains their coefficients, i.e. the classical space $\mathfrak T(G)$ of trigonometric polynomials \cite[27.7]{HR}.

Let $\chi$ be an involutive character of the algebra $\K(G)$. Its restriction onto $\mathfrak T(G)$ is also an involutive character, and it is known \cite[30.5]{HR} that then $\chi(f)=f(t)$ for some point $t\in G$ and all $f\in\mathfrak T(G)$. As $\mathfrak T(G)$ is dense in $\K(G)$ by lemma \ref{K(G)-properties}(5), this is true for all $f\in\K(G)$ as well.
\epr

\blm\label{K(N)-spectrum}
In the notations \eqref{SIN-extension}, every nonzero involutive continuous character $\chi$ of the algebra $\K(N)$ is of the type $\chi(f)=f(t)$ for some point $t\in N$.
\elm
\bpr
By lemma \ref{K(N)-tensor}, the algebras $\K(V)$ and $\K(K)$ may be considered as subalgebras in $\K(N)$ (since each of them contains the constant 1 function). Restrictions of $\phi$ onto $\K(V)$ and $\K(K)$ are characters of these algebras. By lemmas \ref{C(G)-abel}, \ref{K(G)-compact} there are $v_0\in V$, $\kappa_0\in K$ such that $\phi(f)=f(v_0)$ for all $f\in\K(V)$, $\phi(f)=f(\kappa_0)$ for all $f\in\K(K)$. Consequently, $\phi(f)=f(v_0\kappa_0)$ for all $f$ in $\K(V)\cdot\K(K)$, and by lemma \ref{K(N)-tensor}, since $\phi$ is continuous, for all $f\in\K(G)$ as well.
\epr

\btm\label{K(G)-spectrum}
Let $G$ be a Moore group, and let $H$ be the set of nonzero involutive continuous characters of the algebra $\K(G)$. Then every $h\in H$ is of the type $h(f)=f(t)$ for some point $t\in G$. If $H$ is endowed with the weak topology generated by $\K(G)$, then this identification is a homeomorphism.
\etm
\bpr
By definition, $H$ is subset of $\K(G)^\star=\hatC(G)^{\star\star}$, but since $\hatC(G)$ is stereotype, we have $H\subset \hatC(G)$.

Consider the canonical projection $p:G\to D$ (in notations \eqref{SIN-extension}). It generates a homomorphism $p:M_c(G)\to M_c(D)$, which is by universality extended to $p: \hatC(G)\to \hatC(D)$ (we will denote all these maps by the same letter). The conjugate map $p^*:\mathcal K(D)\to\mathcal K(G)$ acts on $f\in \K(D)$ as $(p^* f)(t) = f(p\,t) = f(tN)$, for any $t\in G$. In particular, for the indicator function $I_s$ of a point $s\in D$ we have $p^*I_s=I_{sN}$. One can notice that this is a homomorphism (with pointwise multiplicaton); thus, $ph $ is a character of $\mathcal K(D)$ for every $h\in H$.

For the discrete group $D$, $\hatC(D)=C^*(D)$ and $\mathcal K(D)$ is the Fourier-Stieltjes algebra $B(D)$ (see lemma \ref{discrete}; the latter algebras equal as sets but may have different topologies). The class of Moore groups is closed under passing to the quotients \cite[Vol.~II, p.~1452]{palmer}, so $D$ is a Moore group too. It is known that every discrete Moore group is a finite extension of an abelian group \cite[Theorem 12.4.26 and p.~1397]{palmer}; in particular, it is amenable. It follows that $B(G)=A(G)$, the Fourier algebra of $D$.

By Eymard's theorem, every character of $A(D)$ is evaluation at a point $d\in D$, so we have $ph =\delta_d$ for some $d\in D$. Now for any $s\in D$ we have $h(I_{sN})= h (p^* I_s)= (ph )(I_s)=\delta_d(I_s)$, so $h (I_{sN})=1$ if $sN=dN$ and 0 otherwise.

For any $f\in\mathcal K(G)$, we have then $h (f) = h (f)h (I_{dN})=h (fI_{dN})$. For $t\in G$, $f\in\mathcal K(G)$, let us denote by $\trans tf$ the function $\trans tf(x)=f(xt^{-1})$. By lemma \ref{trans_K(N)}, $I_{dN}\mathcal K(G) = \trans d{\mathcal K(N)}\equiv\{\trans df: f\in \mathcal K(N)\}$.

Then $\trans {d\,}{(\!fg)}=\trans df\,\trans dg$, so $\psi(f)=h (\trans {d}f)$ is a character of $\K(N)$. By lemma \ref{K(N)-spectrum}, $\psi(f)=f(\nu)$ for some $\nu\in N$. It follows that $h (f)=\psi(\trans {d^{-1}\!}f) = (\trans {d^{-1}\!}f)(\nu) = f(\nu d)$ for $f\in\trans d{}\K(N)$, but for arbitrary $f\in\mathcal K(G)$ we have also $h (f) =h (fI_{dN}) = (fI_{dN})(\nu d) = f(\nu d)$, and this proves the equality $G=H$.

To show that this is homeomorphism when $H$ is endowed with the weak topology of $\hatC(G)$, it is equivalent to prove that $\K(G)$ is a regular algebra of functions on $G$, i.e. it separates points and closed sets. First of all, $\K(G)$ contains $I_{tN}$ for any $t\in G$, what implies that every coset $tN$ is open in $H$. So we can further restrict ourselves to $G=N=\R^n\times K$. On the compact subgroup $K$, the imbedding into $H$ is continuous and injective, so this is a homeomorphism, and it implies that $\K(K)$ is a regular algebra. For $\R^n$, this can be verified directly since $\K(\R^n)$ is the Fourier transform of $M_c(\R^n)$; it contains, for example, $f(t)=(\exp(i\e (t_1+\dots t_n))-1)/(i^n t_1\cdots t_n)$ for any $\e>0$ (the Fourier transform of the indicator of the cube $[0,\e]^n$). These functions separate 0 from any closed set in $\R^n$. Now, since on $N$ we have the product topology, the statement follows from the fact that $\K(\R^n)\cdot\K(K)\subset \K(N)$.
\epr

We can now return to the diagram \eqref{quadrat-intro}:
\btm\label{SIN-reflexive}
For any Moore group $G$, $\widehat{\hatC(G)}=C(G)$, and the dual structure of a $\odot$-Hopf *-algebra coincides with the natural structure of $C(G)$. This is illustrated by the following diagram:
$$
 \xymatrix @R=1.pc @C=2.pc
 {
 C(G)\ar@{|->}[rrr]^{\star}
 &
 & &
 M_c(G) \ar@{->}[dd]^{C^*\rm -env}
 \\
 & & &
 \\
 \K(G) \ar@{->}[uu]^{C^*\rm -env}
 & & &
 \hatC(G)\ar@{|->}[lll]_{\star}
 }
$$
\etm
\bpr
Let $j$ denote the canonical inclusion of $G$ into $\hatC(G)$. This map is continuous but it need not be a homeomorphism. If $K\subset G$ is compact, then $j(K)$ is compact in $\hatC(G)$; the seminorm $p_K(f)=\sup_{\phi\in j(K)} |\phi(f)|=\sup_{\phi\in K} |f(t)|$ on $\K(G)$ has the \cst-property and is continuous by definition of topology on $\K(G)$.

Conversely, let $p$ be a continuous \cst-seminorm on $\K(G)$. Then $C_p^*(\K(G))$ is a commutative \cst-algebra, thus isomorphic to $C(K)$, where $K$ is its spectrum, i.e. the set of its involutive continuous characters in the weak-* topology. If $i_p$ is the canonical map from $\K(G)$ to $C_p^*(\K(G))$, then $i_p^*$ maps continuously $K$ to $G$, which is the spectrum of $\K(G)$ by theorem \ref{K(G)-spectrum}. Moreover, $i_p^*$ is injective since $i_p(\K(G))$ is dense in $C_p^*(\K(G))$. By a well-known theorem, this is a homeomorphism of $K$ onto $i_p^*(K)$. Thus, $\K(G)^\Diamond$ and $C(G)$ are inverse limits of equal systems of algebras:
$$
\K(G)^\Diamond = \liminv_{p\in\mathcal P(\K(G))} C_p^*(\K(G)) = \liminv_{K\subset G} C(K) = C(G).
$$
The Hopf structure of $C(G)$ extends that of $\K(G)$, what can be verified directly by evaluating on $\delta$-functions.
\epr

Now, from the previous results follow:
\btm\label{duality_on_Moore_category}
There is a category $\mathfrak{Moore}$ which objects are $\odot$-Hopf *-algebras reflexive with respect to $\widehat{\;}$, and the morphisms are morphisms of $\odot$-Hopf *-algebras. $\mathfrak {Moore}$ contains the algebras $C(G)$ and $\hatC(G)$ for all Moore groups $G$. The mapping $A\mapsto \widehat A$ is a contravariant duality functor on $\mathfrak{Moore}$.
\etm

\section{More properties of the maximal tensor product}\label{section-max}

In this section, we show that among \cst-algebras, every $\omax$-algebra is also a $\odot$-algebra; in particular,  $\hatC(G)$ for a discrete group is a $\omax$-Hopf algebra if and only if $G$ is a Moore group. This means that we could not widen the class of groups in category $\mathfrak {Moore}$ if we considered $\omax$ instead of $\odot$.

\btm\label{mult_bounded_degree}
A \cst-algebra $A$ is strict if and only if it is a $\omax$-algebra.
\etm

This theorem is proved exactly in the same way as it is done by Aristov \cite{aristov} for minimal tensor products. For completeness, we repeat the argument below, including proof of the following proposition which is stated as obvious in \cite{aristov}.

\bprop\label{a_i}
Let $A$ be a \cst-algebra and let $a_1,\ldots, a_m\in A$ be such that $a_i a_j^* = a_i^* a_j = 0$ if $i\ne j$.
Then
$$
\textstyle\| \sum\limits_{i=1}^m a_i \| \le \max\limits_{i=1,\dots,m} \|a_i\|.
$$
\eprop
\bpr
For any $a\in A$, define $S_0(a)=a$ and $S_n(a) = S_{n-1}(a)^* S_{n-1}(a)$ for $n=1,2,\dots$. Then, by induction, $S(a_i)$ satisfy the same conditions as  $a_i$:
$$
S_n(a_i) S_n(a_j)^* = S_{n-1}(a_i)^* S_{n-1}(a_i) S_{n-1}(a_j)^* S_{n-1}(a_j) = 0 \text{ if } i\ne j,
$$
and similarly $S_n(a_i)^* S_n(a_j) = 0$ if $i\ne j$.

These equalities imply that
$$
\textstyle
\|\sum_i S_n(a_i) \|^2 = \| \sum_{i,j} S_n(a_j)^* S_n(a_i) \| = \| \sum_i S_n(a_i)^*S_n(a_i)\| = \| \sum_i S_{n+1}(a_i)\|,
$$
so that by induction
$$
\textstyle\|\sum a_i\|^{2^n} =  \| \sum S_n(a_i)\|.
$$
Now, as $\|S_n(a)\|=\|S_{n-1}(a)\|^2$, we have $\|S_n(a)\|=\|a\|^{2^n}$ and
$$
\textstyle\|\sum a_i\| \le  \left( \sum \|S_n(a_i)\| \right)^{1/2^n} = \left( \sum \|a_i\|^{2^n} \right)^{1/2^n}.
$$
The right-hand side tends to $\max \|a_i\|$ as $n\to\infty$, what proves the proposition.
\epr
\bpr[Proof of the theorem \ref{mult_bounded_degree}.]
If $A$ is strict, then its multiplication $\m$ extends to a continuous operator from $A\underset{\min}\otimes A$ to $A$;
since there is a continuous map from $A\omax A$ to $A\underset{\min}\otimes A$, the multiplication is continuous on
$A\omax A$ also.

For the inverse statement, we follow reasoning of Aristov. In \cite{aristov-sm} it is proved that a \cst-algebra
is strict if and only if it is of bounded degree. Thus, we need to show that every $\omax$-algebra is of bounded degree.
Suppose that $A$ is a $\omax$-algebra and there is an irreducible representation $\pi$ of $A$ of dimension at least $n$.
Let $\{e_i\}_{i=1}^n$ be an orthonormal system in the space $H$ where $\pi$ acts, and let $\e_i$ be the corresponding coordinate functionals: $\e_i(v)=\langle v,e_i\rangle$, $v\in H$. The operator $q=\sum\limits_{i=1}^n \e_i e_i$ is the orthogonal projector onto the span of $e_1,\dots,e_n$, and in particular it is self-adjoint.

By \cite[Lemma 1]{aristov} there exist elements $a_i\in A$, $i=1,\ldots,n$ (in notations of \cite{aristov},
$a_i=v_{i1}$) such that
\begin{gather*}
\|a_i\|=1; \\
a_i a_j^* = a_i^* a_j =0 {\rm \ if\ } i\ne j; \\
q\pi(a_i) q = \pi(a_i)q = \e_1 e_i.
\end{gather*}
Take $x=\sum_i a_i^*\otimes a_i\in A\omax A$. Then
\begin{align*}
q\pi(\m(x))q &= \sum_i q\pi(a_i^*a_i)q = \sum_i \big( \pi(a_i)q\big)^* \pi(a_i)q
 = \sum_i \big( \e_1 e_i \big)^* \e_1 e_i   \\
 &= \sum_i \e_i e_1 \circ \e_1 e_i
 = \sum_i \e_1 e_1 = n \e_1 e_1.
\end{align*}
Thus, $n = \|q\pi(\m(x))q\| \le \|\m(x)\|$. At the same time, the elements $a_i^*\otimes a_i$ satisfy conditions of proposition \ref{a_i}:
$$
(a_i^*\otimes a_i)^* (a_j^*\otimes a_j) = a_i a_j^*\otimes a_i^* a_j,
\quad (a_i^*\otimes a_i) (a_j^*\otimes a_j)^* = a_i^* a_j\otimes a_i a_j^*,
$$
and this is zero is $i\ne j$. By proposition \ref{a_i}, $\|x\|\le\max \|a_i^*\otimes a_i\|=1$. This means that $\|\m\|\ge n$, so if $\m$ is a continuous operator then
$A$ must be of bounded degree.
\epr

Similarly to lemma \ref{nomult}, we have a corollary:
\bcor
For a discrete group $G$ which is not a Moore group, $\hatC(G)$ is not an $\omax$-algebra.
\ecor
\bpr
By lemma \ref{discrete}, $\hatC(G)=C^*(G)$ is a \cst-algebra. If it is a $\omax$-algebra, then by theorem \ref{mult_bounded_degree} it is strict, and in particular, all its irreducible representations are finite-dimensional; this means that $G$ is a Moore group.
\epr

{\bf Acknowledgements.} I would like to thank
S. Akbarov, O. Aristov, V. M. Manuilov, A. Pirkovskii, T. Shulman and A. Van Daele
for numerous discussions and advice. I also thank the referee for valuable comments.

\end{document}